\documentclass{article}

\input{tcilatex}
\begin{document}

\begin{center}
\textbf{Comment on "Sequential Monte Carlo for Bayesian Computation" (P. Del
Moral, A. Doucet, A. Jasra, \textit{Eighth Valencia International Meeting on
Bayesian Statistics})}

David R. Bickel (Johnston, Iowa, USA)

June 21, 2006
\end{center}

The authors added promising innovations in sequential Monte Carlo
methodology to the arsenal of the Bayesian community. Most notably, their
backward-kernel framework obviates the evaluation of the importance density
function, enabling greater flexibility in the choice of algorithms. They
also set their work on posterior inference in a more general context by
citing results of the observation that particle filters approximate the path
integrals studied in theoretical physics (Del Moral 2004).

My first question concerns another recent advance in SMC, the use of the
mixture transition kernel%
\[
\overline{K}_{n}\left( \mathbf{x}_{n-1},\mathbf{x}_{n}\right)
=\dsum\limits_{m=1}^{M}\overline{\alpha }_{n,m}\kappa _{m}\left( \mathbf{x}%
_{n-1},\mathbf{x}_{n}\right) ,
\]%
\newline
where $\overline{\alpha }_{n,m}\ $equals the sum of normalized weights over
all particle values that were drawn from the $m$th mixture component at time 
$n-1,$ and $\kappa _{m}\left( \mathbf{x}_{n-1},\mathbf{x}_{n}\right) $ is an
element of the set of $M$ predetermined mixture components (Douc et al. 2006)%
$.$ For example, if the possible transition kernels correspond to
Metropolis-Hastings random walk kernels of $M$ different scales chosen by
the statistician, then the mixture automatically adapts to the ones most
appropriate for the target distribution (Douc et al. 2006). Is there a class
of static inference problems for which the backward-kernel approach is
better suited, or is it too early to predict which method may perform better
in a particular situation?

In his discussion, Carlo Berzuini suggested some opportunities for further
SMC research. What areas of mathematical and applied work seem most
worthwhile?

I thank the authors for their highly interesting and informative paper.

\bigskip 

\textbf{Reference in Invited Paper}

Del Moral, P. (2004) \textit{Feynman-Kac Formulae: Genealogical and
Interacting Particle Systems with Applications}, Series: Probability and its
Applications, New York: Springer.

\textbf{Additional Reference}

Douc, R., Guillin, A., Marin, J.M., and Robert, C.P. (2006) Convergence of
adaptive sampling schemes. \textit{Annals of Statistics} (to appear). 

\end{document}